\documentclass[a4paper]{eccomas_2022}
\usepackage{graphicx}
\usepackage{amsmath}
\usepackage[hidelinks]{hyperref}
\usepackage{cleveref}
\usepackage{amsfonts}
\usepackage{subcaption} 
\usepackage{mathtools} 
\usepackage{soul}

\title{MODEL ORDER REDUCTION OF SOLIDIFICATION PROBLEMS\break ECCOMAS CONGRESS 2022}

\author{Florian Arbes$^{1,2}$, Øyvind Jensen$^1$, Kent-Andre Mardal$^2$ and Jørgen S. Dokken$^3$}

\heading{Florian Arbes, Øyvind Jensen, Kent-Andre Mardal and Jørgen S. Dokken}

\address{$^{1}$ Institute for Energy Technology (IFE)\\
Instituttveien 18 2007 Kjeller, Norway\\
e-mail: \{florian.arbes, oyvind.jensen\}@ife.no
\and
$^{2}$ Department of Mathematics, University of Oslo (UiO)\\
Moltke Moes vei 35, 0851 Oslo, Norway\\
email: kent-and@math.uio.no
\and
$^{3}$ Department of Engineering, University of Cambridge\\
Trumpington St, Cambridge CB2 1PZ, United Kingdom\\
email: jsd55@cam.ac.uk
}


\abstract{
Advection driven problems are known to be difficult to model with a reduced basis because of a slow decay of the Kolmogorov $N$-width.
This paper investigates how this challenge transfers to the context of solidification problems and tries to answer when and to what extend reduced order models (ROMs) work for solidification problems.
In solidification problems, the challenge is not the advection per se, but rather a moving solidification front.
This paper studies reduced spaces for 1D step functions that move in time, which can either be seen as advection of a quantity or as a moving solidification front.
Furthermore, the reduced space of a 2D solidification test case is compared with the reduced space of an alloy solidification featuring a mushy zone.
The results show that not only the PDE itself, but the smoothness of the solution is crucial for the decay of the singular values and thus the quality of a reduced space representation.
}

\begin{document}

\section{INTRODUCTION}
Solidification processes in the material processing industry are complex and need to be modelled numerically.
The underlying multi-physics equations for fluid flow, heat transfer, and phase change must be coupled to one another \cite{beckermann_mathematical_1993}.
After discretizing in time, the underlying equations are often split up to speed up the numerical computation of the solution.
Velocity components and the pressure are solved for in a segregated manner.
There is a number of so called fractional step algorithms to solve the incompressible Navier-Stokes equations, see e.g. \cite{langtangen_numerical_2002,turek_efficient_2014} for overviews.
Often, three steps are involved \cite{mortensen_oasis_2015}: firstly, one solves for a tentative velocity vector $\textbf{u}^I$ component-wise based on a pressure guess $p^*$. Secondly, a pressure correction $\varphi = p^{n-1/2}-p^*$ is computed and lastly a velocity correction is computed.
The $k$-th velocity component is denoted with $u^I_k$.
Step one and two may be repeated iteratively to achieve higher accuracy.
Simo and Armero \cite{simo_unconditional_1994} outline a general fractional step algorithm as followed:

\begin{align}
\frac{u^I_k-u^{n-1}_k}{\Delta t} + B^{n-1/2}_k &= \nu \nabla ^2 \widetilde{u}_k - \nabla_k p^* + f^{n-1/2}_k \text{ for } k = 1, ..., d, \label{equ:step1}\\
\nabla ^2 \varphi &= -\frac{1}{\Delta t} \nabla \cdot \textbf{u}^I,
 \label{equ:step2}\\
\frac{u^n_k-u^I_k}{\Delta t} &= - \nabla_k \varphi \text{ for } k = 1, ..., d,
\label{equ:step3}
\end{align}

where $ B^{n-1/2}_k$ denotes the nonlinear convection term evaluated at the midpoint $n-1/2$ of the time interval $\Delta t$ and $d$ the dimensionality of the problem.

\begin{figure}[!ht]
 \centering
 \includegraphics[width=8cm]{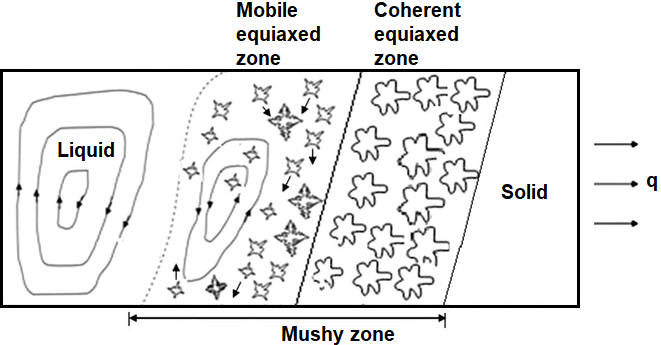}
\caption{The mushy zone is a mixture of solid and liquid that leads to a reduced fluid flow (adapted from \cite{kumar_grain_2011}).}
 \label{fig:mushy-zone}
\end{figure}

For alloy solidifications, the interface between solid and liquid is is not sharp.
The different components in an alloy composition have different melting points, therefore mobile equiaxed grains appear in the melt when cooling down.
While solidifying, these grains form a coherent structure before the composition fully solidifies at a lower temperature (c.f. \cref{fig:mushy-zone}).
The so called mushy zone is defined by a temperature dependent mix of solid and liquid.
Non-linearities like that make for challenging numerical modelling and high computational costs, which is problematic when either fast predictions are needed or when an application is studied for varying parameters.

To address computational challenges, reduced order models (ROMs) have been developed.
ROMs aim to approximate solutions efficiently at little or no expense of accuracy.
Fast and reliable predictions open up new possibilities to study a numerical problem for varying parameters.
Reduced order modelling refers to replacing a high fidelity problem with a lower complexity numerical model \cite{quarteroni_reduced_2015}.
PDEs are solved by discretizing the domain, which then leads to high-dimensional numerical models with millions of degrees of freedom \cite{richardson_scalable_2019}.
The solutions resulting from different values of a parameter often suggest an underlying structure that can be described by a model with fewer degrees of freedom.
A ROM aims to evaluate the solution of any set of new parameters at a computational cost that does not depend on the number of degrees of freedom of the high fidelity model.

ROMs are built in what is referred to as an offline stage and employed in an online stage.
In the offline stage, a reduced basis is generated for example with a greedy approach or via a proper orthogonal decomposition (POD).
In any case, the high fidelity full order model (FOM) is evaluated several times for different parameters $\mu \in \mathbb{R}^d $ that characterize the problem.
In case of the solidification problem, the parameters might be cooling rates, wall temperatures, or material properties \cite{zimmerman_monolithic_2018,torabi_rad_theory-training_2020}.
In order to build the ROM using a POD, the discrete solutions $u_{\delta}(\mu)$ (snapshots) need to be collected in a matrix $X$.
Each column in the snapshot matrix corresponds to one snapshot, which means the column index $j$ maps the parameter space $\mu$ to the corresponding snapshots.
Each row in the snapshot matrix corresponds to a node and its associated quantity (velocity in x, y, z, pressure, solid fraction, etc.).
The POD can then be represented by a singular value decomposition (SVD) of the the snapshot matrix $X$ \cite{brunton_data-driven_2019}:

\begin{equation}
\label{SVD}
X = U \Sigma V^T.
\end{equation}

Once the reduced basis is found in the offline phase, the model can be employed in the online phase, which can be done in an intrusive or non-intrusive manner.
For a new parameter, the reduced basis solution $v_{rb}(\mu)$ can be found by projecting the governing equations onto the reduced space, also known as a Galerkin projection.
A Galerkin projection guarantees that the error of the reduced space approximation is orthogonal to the reduced space, meaning the result is the best possible solution in the reduced space.
However, this approach often suffers from stability issues and runs into efficiency problems when it comes to non-linear problems \cite{chen_blackbox_2012,xiao_non-intrusive_2016}.
There are methods that try to overcome these limitations, the interested reader is referred to books on the subjects \cite{brunton_data-driven_2019,quarteroni_reduced_2015,benner_model_2017,hesthaven_certified_2015}.
Since the "POD-Galerkin" approach depends on the equations of the FOM, existing code has to modified and expanded, which is cumbersome and prone to error.
Multi-physics problems such as a solidification problem rely on extensive code libraries \cite{mortensen_development_2011}.
The equations of the fluid flow have to be coupled with the equations for the heat transfer as well as those of a the phase change, and non-linearities need to be handled carefully~\cite{beckermann_mathematical_1993}.

We seek a non-intrusive ROM, because these models are purely data-driven surrogate models that do not depend on the FOM.
Non-intrusive reduced order models (NIROM) rely on methods that infer the right singular values from the parameters $\mu$ that corresponds to the snapshot matrix.
This can be done with e.g. a radial basis function (RBF) interpolation ~\cite{xiao_non-intrusive_2015-1}.
A straightforward approach would be a standard interpolation on a regular grid, however, this could require a lot of snapshots if the dimensionality of the parameter space is high.
Xiao et al. proposed Smolyak sparse grids in this context to avoid an extensive amount of snapshots~\cite{xiao_non-intrusive_2016,smolyak_quadrature_1963}.
Deep neural networks have successfully been used for the mapping of the right singular values as well as for uncertainty quantification
~\cite{pawar_deep_2019,jacquier_non-intrusive_2021}.

For a reduced order model it is important to know how well the solution manifold $M_{\delta}$ can be approximated.
The solution manifold is the set of all solutions of the problem for any choice of the parameters.
The error arising from a reduced space approximation $v_{rb}$ obtained through a Galerkin projection is measured with the Kolmogorov $N$-width $d_N$~\cite{hesthaven_certified_2015}:

\begin{equation}
d_N(M_{\delta}) = \inf\limits_{\mathbb{V}_{rb}} \sup\limits_{u_{\delta} \in M_{\delta}} \inf\limits_{v_{rb} \in \mathbb{V}_{rb}} ||u_{\delta}(\mu)-v_{rb}(\mu)||_\mathbb{V}.
\label{equ:kolmogorov}
\end{equation}

The error depends on the number of basis functions $N$ that span the subspace, hence the name "$N$-width".
A fast decay of the $N$-width means the underlying problem is low-dimensional and it can be approximated with a fewer basis functions.
The $N$-width is problem-dependant and depends on the coercivity and continuity constants of the bilinear form (c.f. Céa´s Lemma).
Ohlberger \cite{ohlberger_reduced_2016-2} points out that advection driven phenomena suffer from a very slow decay of the Kolmogorov $N$-width.
It will be demonstrated, that discontinuities or jumps of a quantity in space that move in time correlate with a slow decay of the singular values.
Jumps appear in the solution snapshots for solidification problems too: 
Due to the abrupt variation in material properties associated with the phase change, the velocity field  $\mathbf{u}$ may vary strongly over a short distance
These non-stationary jumps in the snapshots do not propagate by advection, but rather by the complex physics of the solidification process.
However, the POD is ignorant of how the non-stationary jumps originally came into the snapshots and we find it fruitful to analyse solidification problems using results developed for advection problems.
This paper aims to investigate how well reduced spaces can approximate solidification problems and how to extend reduced order models for such problems.

\section{1D TEST CASES}

It is widely known, that pure Poisson problems have an underlying structure with a low rank that can be exploited using a ROM that is built on a linear subspace found with a POD \cite{ohlberger_reduced_2016-2,brunton_data-driven_2019}.
The Kolmogorov $N$-width decays exponentially fast.
In contrast to that, advection driven problems often suffer from a slow decay of the Kolmogorov $N$-width.

\begin{figure}
\centering
\begin{subfigure}{.5\textwidth}
 \centering
 \includegraphics[width=1\linewidth]{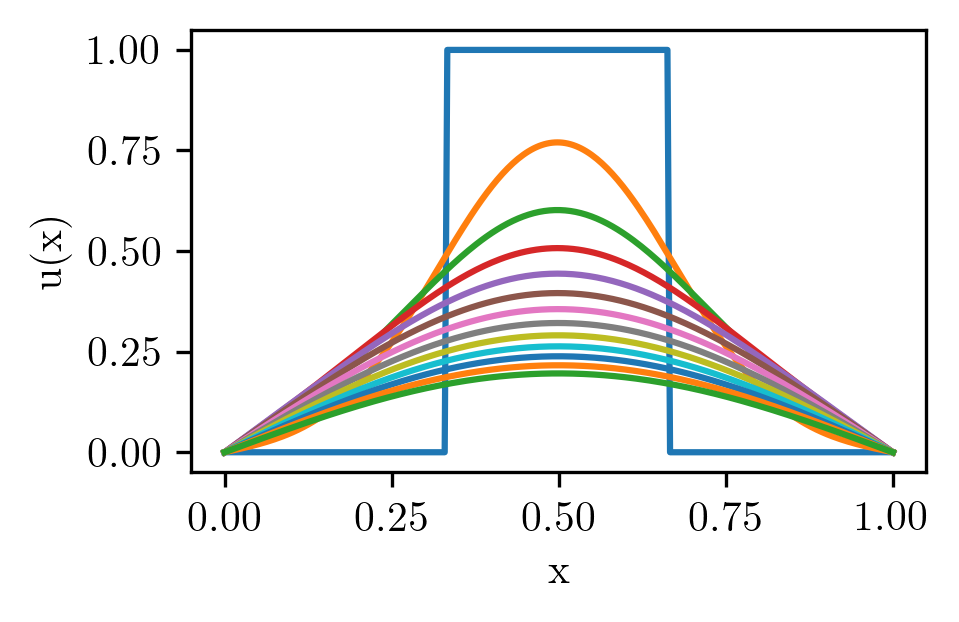}
 \caption{Solutions to the 1D heat equation}
 \label{fig:poisson_snapshots}
\end{subfigure}%
\begin{subfigure}{.5\textwidth}
 \centering
 \includegraphics[width=1\linewidth]{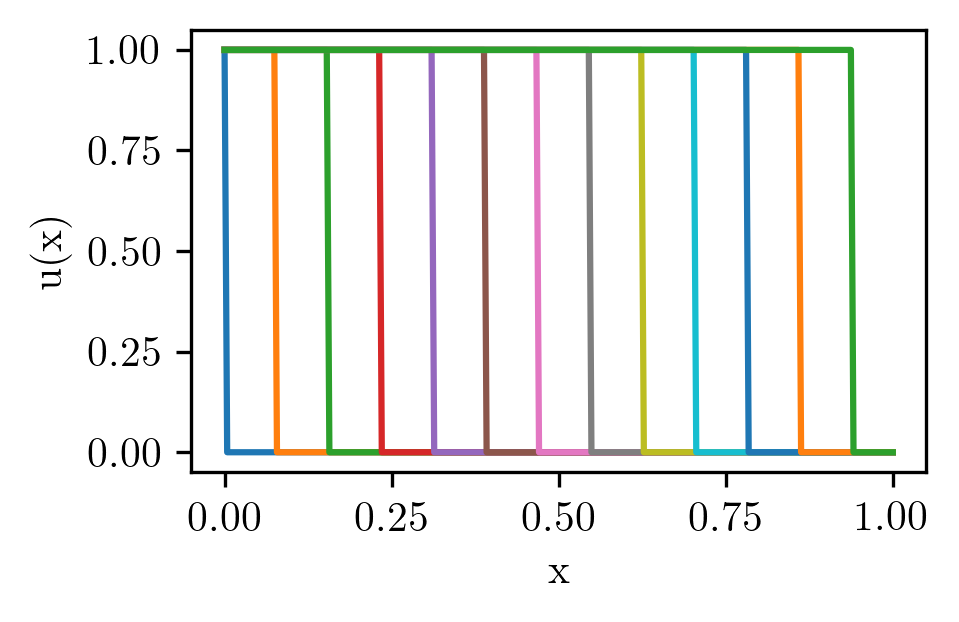}
 \caption{Advected jump discontinuity}
 \label{fig:advection1_snapshots}
\end{subfigure}%

\caption{Selected snapshots of two classical time dependant problems.}
\label{fig:poisson_advection_SS}
\end{figure}

In order to give a demonstration and a reference, two simple time dependant 1D problems have been set up.
The first solves the heat equation, with a rectangular function as initial condition that fulfills homogeneous Dirichlet boundary conditions:
\begin{align}
\begin{split}
\frac{\partial u}{\partial t} - \nabla \cdot \left ( \alpha \nabla u \right ) &= 0 \quad {\rm in} \ \Omega,\label{equ:poisson} \\
u &= 0 \quad {\rm on} \ \partial \Omega.
\end{split}
\end{align}
The thermal conductivity $\alpha$ as well as the timestep $\Delta t$ were chosen to be constant ($\alpha = 1.0$ and $\Delta t = 0.001$).
After discretizing in time using forward Euler, the equation is solved for each timestep using first order Lagrange elements.
The problem is solved using 256 nodes in space and the first 128 snapshots are collected in a snapshot matrix.
Some of the resulting snapshots are shown in~\cref{fig:poisson_snapshots}.

The second reference problem is an advected step function, as described in \cite{ohlberger_reduced_2016-2}, where the solution is given as:
\begin{equation}
 u(x, t) = \begin{dcases*}
 1, & if $ x \leq t $,\\
 0, & \text{otherwise}.
 \end{dcases*}
\label{equ:step}
\end{equation}

In \cref{equ:step}, $x$ is the spatial coordinate of the nodes and $t\in[0, 1]$.
For the sake of consistency, the number of nodes is chosen to be $N=256$ and there are 128 snapshots being collected.
Some of the resulting snapshots are shown in \cref{fig:advection1_snapshots}.

\begin{figure}
\centering
\begin{subfigure}{.48\textwidth}
 \centering
 \includegraphics[width=1\linewidth]{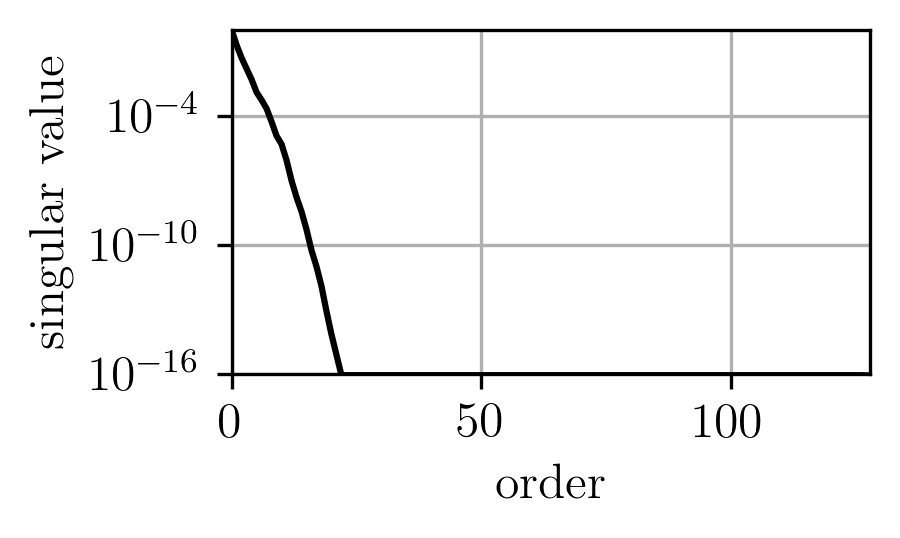}
 \caption{Normalized decay of the singular values of the heat diffusion problem (snapshots shown in \cref{fig:poisson_snapshots})}
 \label{fig:poisson_SV}
\end{subfigure}%
\hfill
\begin{subfigure}{.48\textwidth}
 \centering
 \includegraphics[width=1\linewidth]{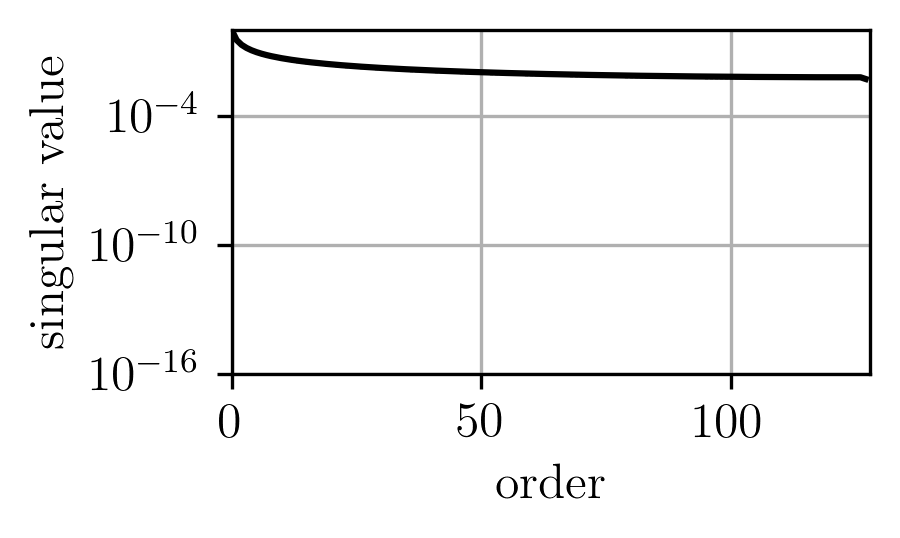}
 \caption{Normalized decay of the singular values of the advected jump (snapshots shown in \cref{fig:advection1_snapshots})}
 \label{fig:advection1_SV}
\end{subfigure}%

\caption{Normalized decay of the singular values of two classical time dependant problems.}
\label{fig:poisson_advection1_SV}
\end{figure}

The resulting snapshot matrices are decomposed using POD.
As expected, the solution snapshots from the heat equation show an exponentially fast decay of the singular values (see \cref{fig:poisson_SV}).
In contrast to that, the singular values of the snapshot matrix containing snapshots from the advection reference case decay slowly (see \cref{fig:advection1_SV}).
This fits nicely to theoretical analysis of the Kolmogorov $N$-width. In fact, Ohlberger and Rave~\cite{ohlberger_reduced_2016-2} have shown, that the Kolmogorov $N$-width decays with $\frac{1}{2\sqrt{n}}$.
Greif and Urban \cite{greif_decay_2019} extended that proof to hyperbolic wave equation with discontinuous initial conditions.

Advected jump discontinuities are rarely seen in real world problems.
In case of no-slip boundary conditions, the velocity profile of a fluid flow problem is not expected to feature a jump discontinuity in its profile.
Even if the boundary layer is thin, the velocity field can be expected to be smooth.
However, the boundary layers might move through the domain, e.g. with a progressing solidification.
Advected quantities such as temperature or enthalpy are also expected to be smooth because of diffusivity.

\begin{figure}
\centering
\begin{subfigure}{.5\textwidth}
 \centering
 \includegraphics[width=1\linewidth]{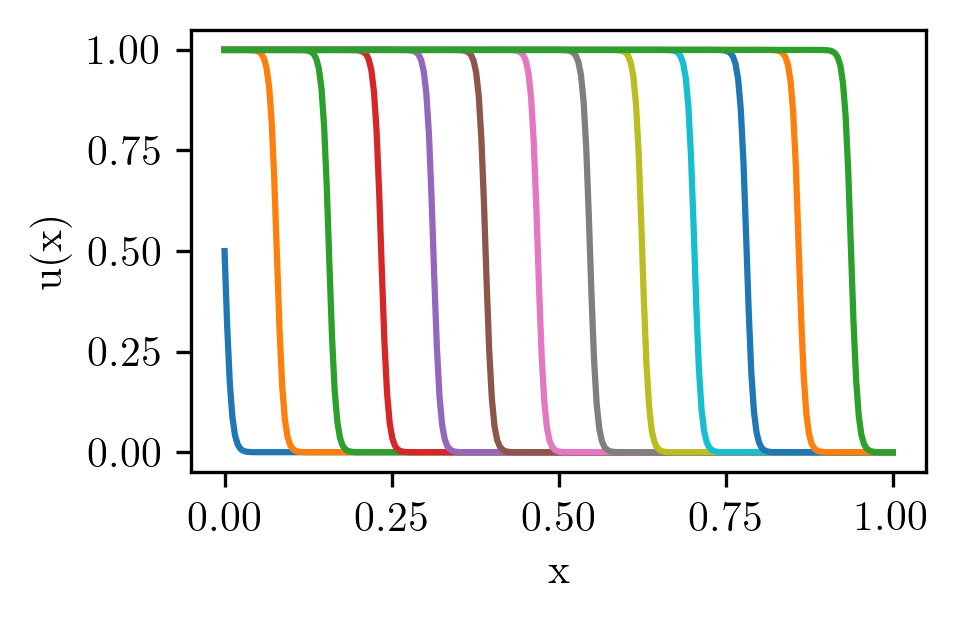}
 \caption{Advected steep sigmoid function}
 \label{fig:advection2_snapshots}
\end{subfigure}%
\hfill
\begin{subfigure}{.5\textwidth}
 \centering
 \includegraphics[width=1\linewidth]{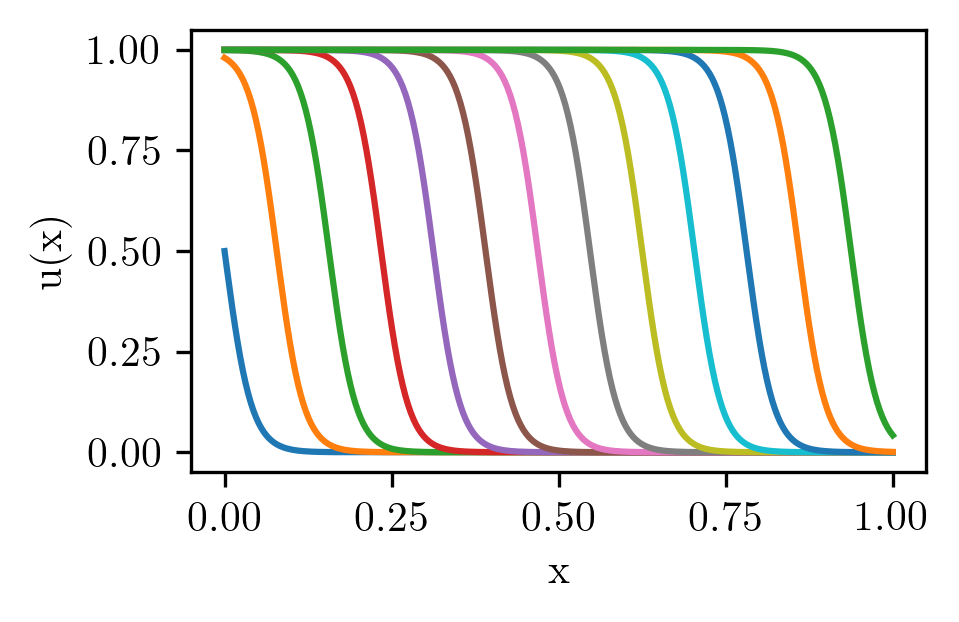}
 \caption{Advected stretched sigmoid function}
 \label{fig:advection3_snapshots}
\end{subfigure}%
\caption{Selected snapshots of two classical time dependent problems.}
\label{fig:advection23_snapshots}
\end{figure}

\begin{figure}
\centering
\begin{subfigure}{.48\textwidth}
 \centering
 \includegraphics[width=1\linewidth]{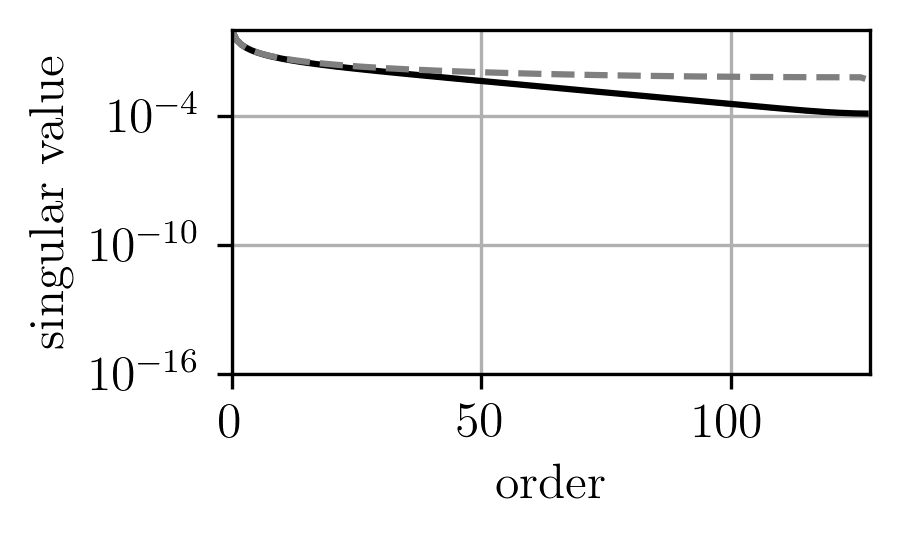}
 \caption{Normalized decay of the singular values of the advected steep sigmoid function (snapshots shown in \cref{fig:advection2_snapshots})}
 \label{fig:advection2_SV}
\end{subfigure}%
\hfill
\begin{subfigure}{.48\textwidth}
 \centering
 \includegraphics[width=1\linewidth]{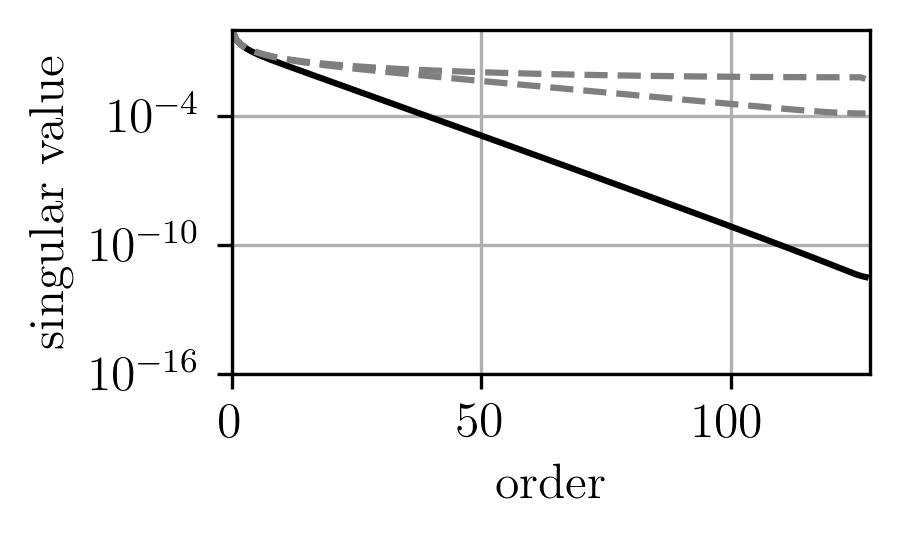}
 \caption{Normalized decay of the singular values of the advected stretched sigmoid function (snapshots shown in \cref{fig:advection3_snapshots})}
 \label{fig:advection3_SV}
\end{subfigure}%
\caption{Normalized decay of the singular values of two advection problems.}
\label{fig:advection23_SV}
\end{figure}

In contrast to the advected jump discontinuity, an advected smooth step function in the shape of a sigmoid function (c.f. \cref{fig:advection2_snapshots}) does seem to decay exponentially fast (c.f. \cref{fig:advection2_SV}).
With decreasing steepness (c.f. \cref{fig:advection3_snapshots}), the decay of the singular values is faster.

We have shown that the decay of the singular values does not only depend on the PDE, but also on the slope and smoothness of the solution.
These results are relevant in the context of solidification of alloys, as the naturally occurring mushy zone has has the effect of smoothing out the jump in the velocity field and thus enhances the decay of the singular values.
The next chapter will discuss how these results transfer to a 2D solidification with or without a mushy zone.

\section{SOLIDIFICATION TEST CASES}
\begin{figure}[!ht]
 \centering
 \includegraphics[width=1.0\textwidth]{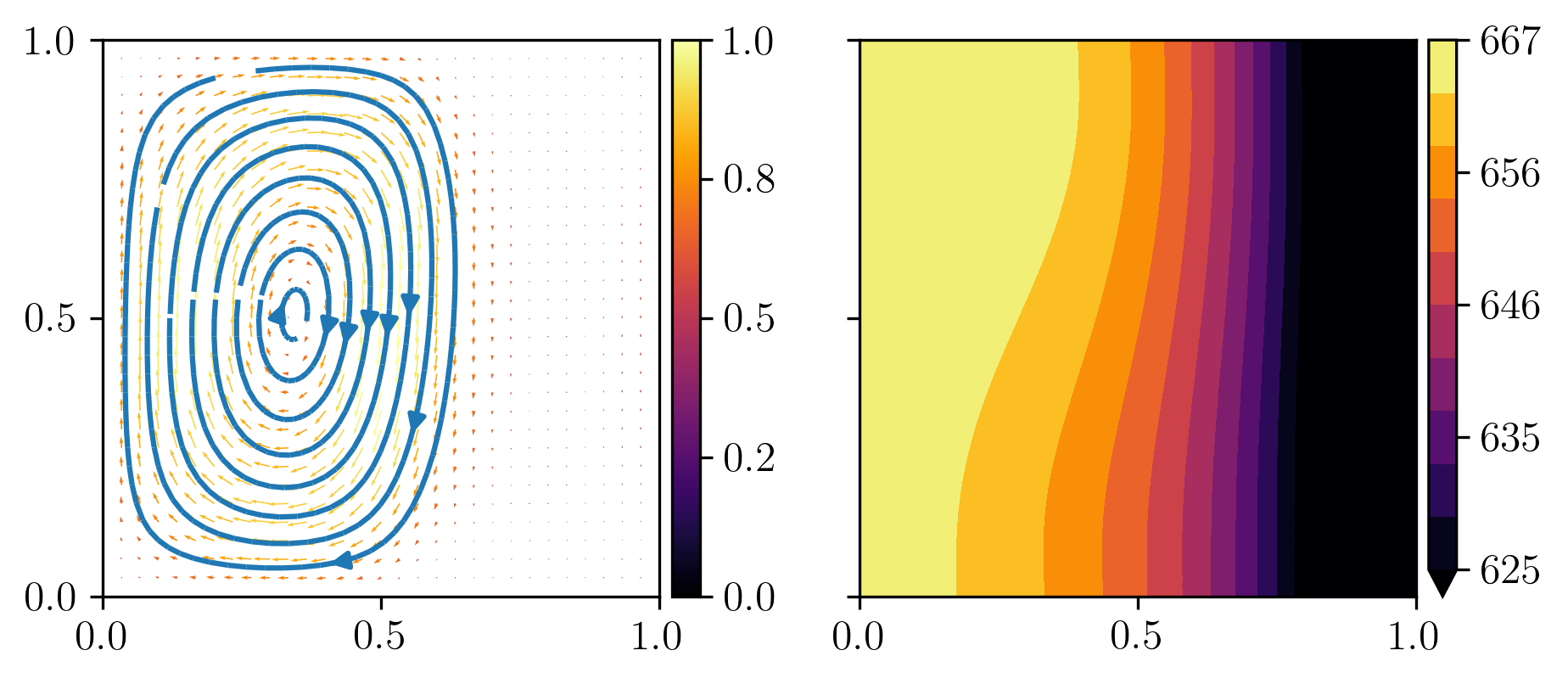}
\caption{Velocity (normalized values, plotted on a coarse grid) and temperature field of a freezing unit square featuring a mushy zone (snapshot in time).}
 \label{fig:freezing_cavity}
\end{figure}

A 2D solidification test case has been set up to test the impact of the mushy zone on the decay of the singular values.
The domain is a unit square, that allows heat transfer through the boundary on the right side, as seen in \cref{fig:freezing_cavity}.
Using the Boussinesq approximation, an implicit pressure correction scheme was implemented to solve the problem.
The scheme treats the convection implicitly by using an Adams-Bashforth projected convecting velocity and Crank-Nicolson for the convected velocity.
Standard Taylor-Hood (continuous P2-P1) elements were used.
In this work we model the restrictions on flow due to the solidification through an increase in viscosity.
In case of an alloy solidification, the viscosity $\mu$ is calculated as $\mu=10(T-650)^2$ below $650^\circ C$, where $T$ is the temperature in $^\circ C$.
This allows for some movement in the mushy zone, because there is a diffuse interface between solid and liquid, as displayed in \cref{fig:mushy-zone}.
In case of a pure metal solidification, i.e. without mushy zone, there is a jump in viscosity over several orders of magnitude at the freezing point ($650^\circ C$) to model a sharp interface between solid and liquid.

\begin{figure}[!ht]
 \centering
 \includegraphics[width=8cm]{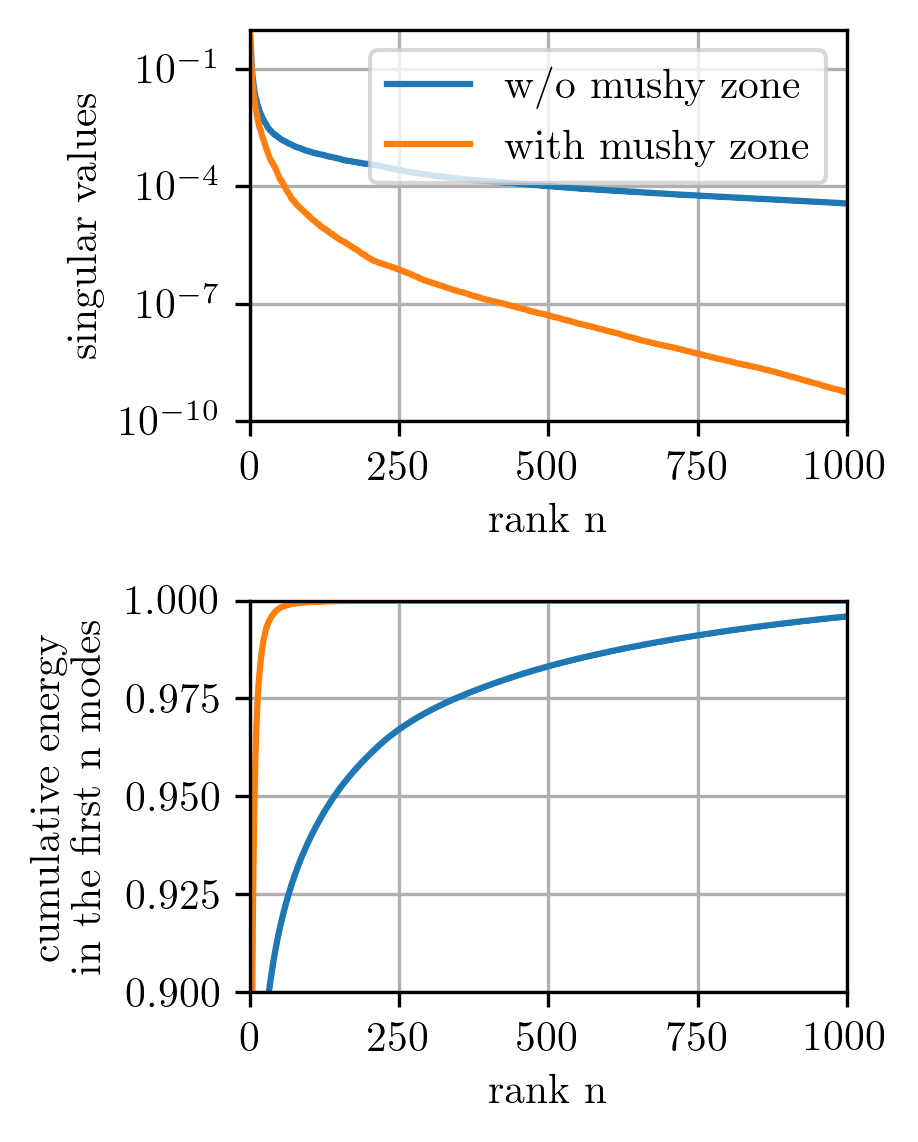}
\caption{Normalized decay of the singular values for a 2D solidification case study with and without a mushy zone.}
 \label{fig:mush_vs_no_mush2}
\end{figure}

For both cases, more than $5 000$ snapshots have been collected.
The subsequent SVD revealed, that the mushy zone has a big impact on the decay of the singular values, c.f. \cref{fig:mush_vs_no_mush2}.
To capture $99.99 \% $ of the energy in $5 000$ snapshots, only $171$ snapshots are needed in case of a mushy zone being present.
If there is no mushy zone, one needs $1527$ snapshots.
These results support the 1D results found in the previous chapter and demonstrate the significance of a mushy zone for a ROM.

\begin{figure}
\centering
\begin{subfigure}{.48\textwidth}
 \centering
 \includegraphics[width=1\linewidth]{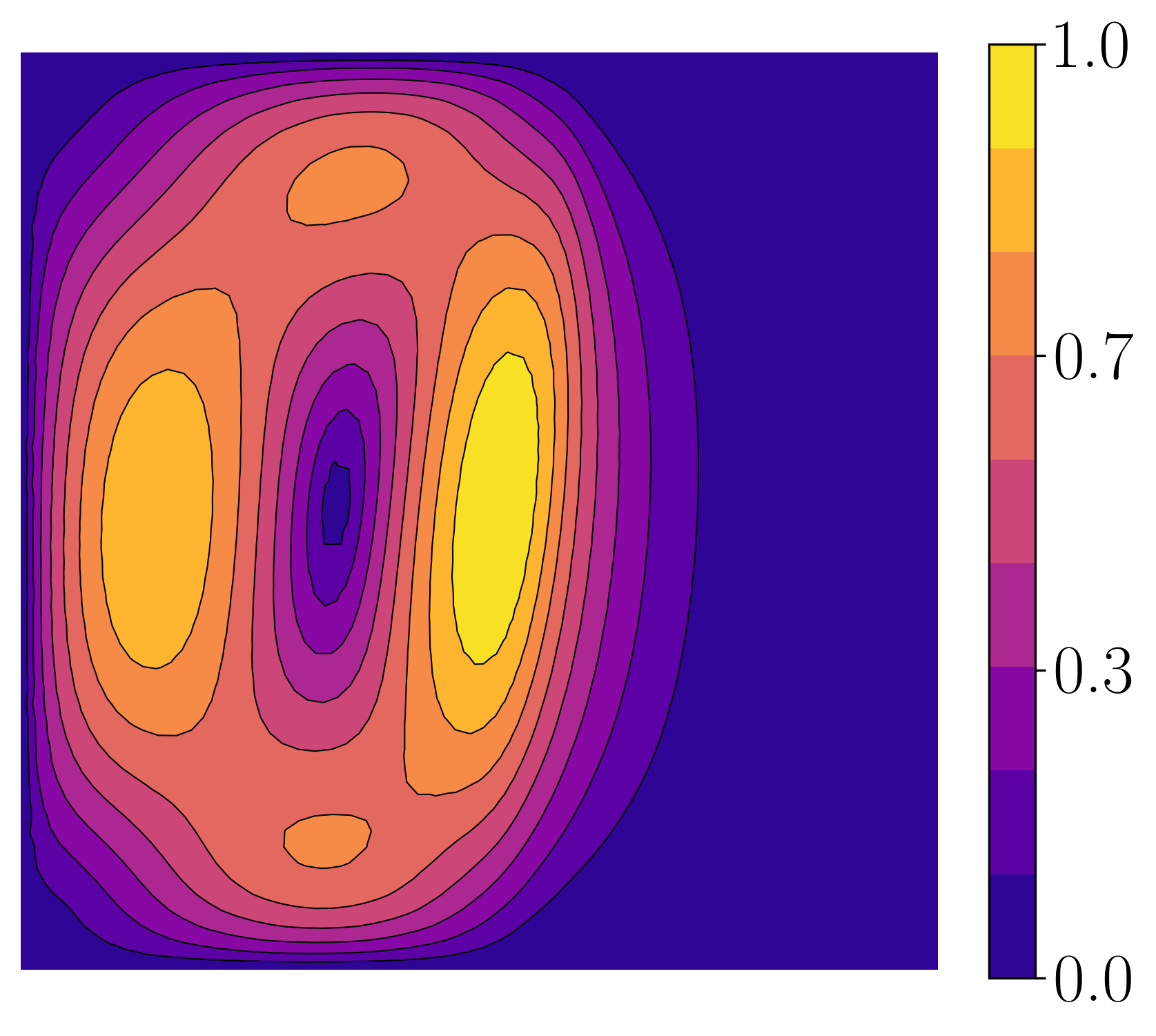}
 \caption{Normalized velocity magnitude of a case featuring a mushy zone (snapshot in time).}
 \label{fig:mushy1}
\end{subfigure}
\hfill
\begin{subfigure}{.48\textwidth}
 \centering
 \includegraphics[width=1\linewidth]{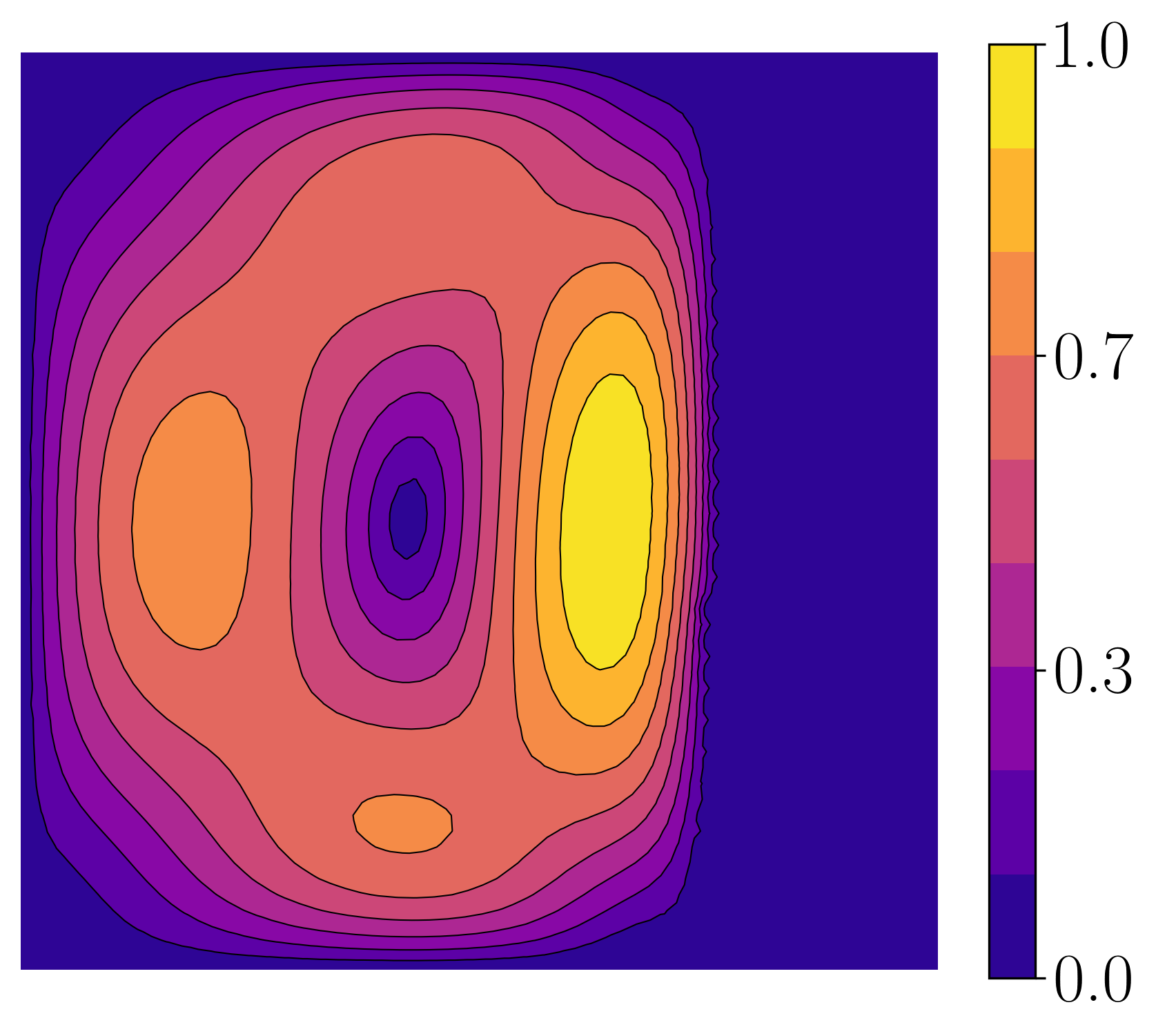}
 \caption{Normalized velocity magnitude of a case without a mushy zone.}
 \label{fig:mushy2}
\end{subfigure}
\caption{Magnitude of the velocity.}
\label{fig:mushy12}
\end{figure}

It is interesting to note, that the velocity field seems to be responsible for the slow decay of the singular values.
\Cref{fig:mushy12} shows a contour plot of the normalized velocity magnitude for both cases.
The temperature driven convection is clearly visible, as well as the solid, that does not allow for any movement.
The mushy zone is also visible in \cref{fig:mushy1}, as the magnitude of the velocity increases gradually and the distance between the solid and the largest velocity is large.
\Cref{fig:mushy1} shows in contrast a smaller proximity of small and large values of the velocity magnitude, because there is no mushy zone.
As seen in the 1D case study, a steep slope in the advected quantity seems to have a negative impact on the decay of the singular values.

\section{ON THE POTENTIAL OF A COMPONENT WISE ROM}
\begin{figure}[!ht]
 \centering
 \includegraphics[width=8cm]{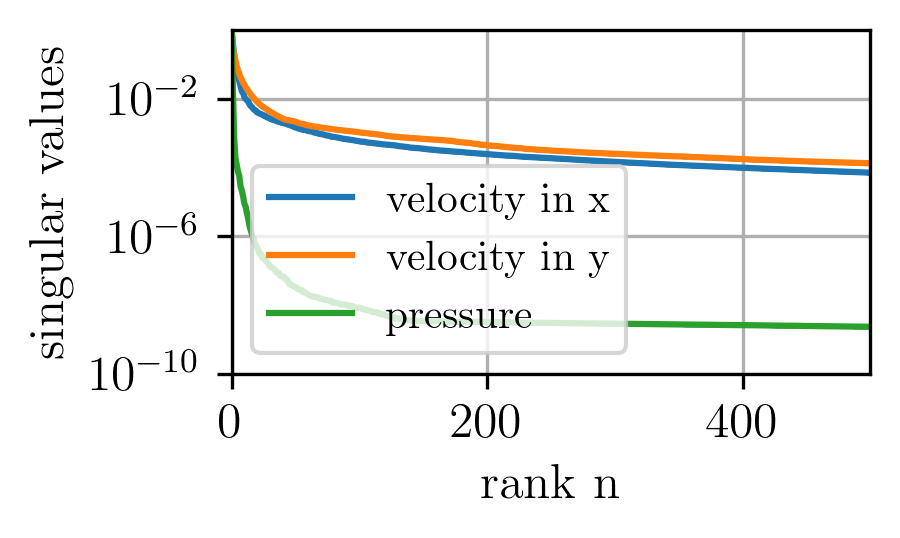}
\caption{Separately computed normalized decay of the singular values for velocity in x- and y-direction and pressure.}
 \label{fig:SV_uvp}
\end{figure}

As seen in the previous chapter, the velocity field of a solidification without a mushy zone shows steep gradients, which leads to a slow decay of the singular values.
Even though velocity, pressure, temperature and viscosity are coupled through the PDE, they are usually computed in a segregated manner.
That motivates to investigate the decay of other quantities separately for the solidification case without a mushy zone.
In \cref{fig:SV_uvp}, the decay of the components is plotted separately.

As expected, singular values for the velocity decays slowly, and the singular values for the pressure decays much quicker.
Since the pressure is computed separately (see \cref{equ:step2}), this opens up the opportunity to build a ROM for the pressure only.
A pressure ROM could potentially be used in the numerical fractional step algorithm of the full order model, because pressure and velocity are treated separately.

\section{CONCLUSIONS}
The quality of metal casts can benefit from fast numerical models that help to prevent defects.
However, the potential for reduced order models is limited, because advection dominated problems suffer from a slow decay of the singular values.
This challenge transfers to metal solidification problems, but not in the same extent to alloy solidification problems.
The numerical experiments conducted suggest that the velocity field moving with the solidification front leads to a slow decay of the singular values.
However, an alloy solidification with a mushy zone allows for much better ROMs since the singular values decay much faster.
The mushy zone allows for some movement in the solid-liquid mixture, which leads to a smoother velocity field.
In this paper, we demonstrated that not only the PDE itself, but the smoothness of the solution determines the decay of the singular values.
This was shown for a simple 1D test case and the results from the 2D solidification test cases suggest that this applies to alloy solidifications too.
Moreover it has been shown, that the pressure decays much faster than the slowly decaying velocity in the metal solidification case.
A pressure ROM opens up opportunities to speed up numerical fractional step algorithms that treat pressure and velocity separately.

\section{ACKNOWLEDGMENT}
The authors gratefully acknowledge the STIPINST funding [318024] from the Research Council of Norway.

\bibliographystyle{biolett}
\bibliography{references}

\end{document}